\documentclass[12pt,leqno]{article}
\usepackage{amsfonts}
\usepackage{amsmath, amsthm, amsfonts, amssymb, color}
\usepackage{mathrsfs}
\usepackage{color}
\theoremstyle{definition}

\newcommand{\scr}[1]{\mathscr #1}
\definecolor{wco}{rgb}{0.5,0.2,0.3}

\numberwithin{equation}{section} \theoremstyle{remark}

\newcommand{\ua}{\uparrow}

\title{{\bf   Ornstein$-$Uhlenbeck Type Processes  on Wasserstein  Spaces}\footnote{Feng-Yu Wang is supported in
 part by  the National Key R\&D Program of China (No. 2022YFA1006000, 2020YFA0712900) and NNSFC (11921001). Panpan Ren is supported by NNSFC(1230018) and  Research Centre for Nonlinear Analysis at The Hong Kong Polytechnic University.}.  }
\author{{\bf  Panpan Ren$^{2)}$ and Feng-Yu Wang$^{1)}$   }\\
\footnotesize{ $^{1)}$ Center for Applied Mathematics, Tianjin University, Tianjin 300072, China}\\
\footnotesize{$^{2)}$ Department of Mathematics,
City University of Hong Kong, Tat Chee Avenue, Hong Kong,  China}\\
\footnotesize{ wangfy@tju.edu.cn, panparen@cityu.edu.hk}}
\begin{document}
\allowdisplaybreaks
\def\R{\mathbb R}  \def\ff{\frac} \def\ss{\sqrt} \def\B{\mathbf
B}
\def\N{\mathbb N} \def\kk{\kappa} \def\m{{\bf m}}
\def\ee{\varepsilon}\def\ddd{D^*}
\def\dd{\delta} \def\DD{\Delta} \def\vv{\varepsilon} \def\rr{\rho}
\def\<{\langle} \def\>{\rangle}
  \def\nn{\nabla} \def\pp{\partial} \def\E{\mathbb E}
\def\d{\text{\rm{d}}} \def\bb{\beta} \def\aa{\alpha} \def\D{\scr D}
  \def\si{\sigma} \def\ess{\text{\rm{ess}}}\def\s{{\bf s}}
\def\beg{\begin} \def\beq{\begin{equation}}  \def\F{\scr F}

\def\e{\text{\rm{e}}} \def\ua{\underline a} \def\OO{\Omega}  \def\oo{\omega}
 \def\tt{\tilde}\def\[{\lfloor} \def\]{\rfloor}
\def\cut{\text{\rm{cut}}} \def\P{\mathbb P} \def\ifn{I_n(f^{\bigotimes n})}
\def\C{\scr C}      \def\aaa{\mathbf{r}}     \def\r{r}
\def\gap{\text{\rm{gap}}} \def\prr{\pi_{{\bf m},\varrho}}  \def\r{\mathbf r}
\def\Z{\mathbb Z} \def\vrr{\varrho} \def\ll{\lambda}
\def\L{\scr L}\def\Tt{\tt} \def\TT{\tt}\def\II{\mathbb I}

\def\M{\mathbb M}\def\Q{\mathbb Q} \def\texto{\text{o}} \def\LL{\Lambda}

\def\to{\rightarrow}\def\l{\ell}\def\iint{\int}\def\gg{\gamma}
\def\EE{\scr E} \def\W{\mathbb W}
\def\A{\scr A} \def\Lip{{\rm Lip}}\def\S{\mathbb S}
\def\BB{\scr B}\def\Ent{{\rm Ent}} \def\i{{\rm i}}\def\itparallel{{\it\parallel}}
\def\g{{\mathbf g}}\def\Sect{{\mathcal Sec}}\def\T{\mathbb T}\def\BB{{\bf B}}
\def\f{\mathbf f} \def\g{\mathbf g}\def\BL{{\bf L}}  \def\BG{{\mathbb G}}
\def\Bd{{D^E}} \def\BdP{D^E_\phi} \def\Bdd{{\bf \dd}} \def\Bs{{\bf s}} \def\GA{\scr A}
\def\Bg{{\bf g}}  \def\Bdd{\psi_B} \def\supp{{\rm supp}}\def\div{{\rm div}}
\def\ddiv{{\rm div}}\def\osc{{\bf osc}}\def\1{{\bf 1}}\def\BD{\mathbb D}\def\GG{\Gamma}
\def\H{\mathbb H}

\maketitle

\begin{abstract}   Let $\scr P_2$  be the space of probability measures on $\R^d$ having finite second moment, and consider   the    Riemannian structure on $\scr P_2$ induced by the intrinsic derivative on the $L^2$-tangent space. By using  stochastic analysis on the tangent space, we construct an Ornstein$-$Uhlenbeck  (OU) type Dirichlet form on $\scr P_2$ whose generator is  formally given by the  intrinsic Laplacian   with a drift.  The log-Sobolev inequality holds and the associated Markov semigroup is $L^2$-compact.
  Perturbations of the OU Dirichlet form are also studied.
\end{abstract}

 \vskip 2cm

 \section{Introduction}

To construct the Ornstein$-$Uhlenbeck process on the Wasserstein space, we introduce an inherent Gauss measure
on the Wasserstein space,  which together with the intrinsic derivative provides an Ornstein$-$Uhlenbeck type Dirichlet form. Up to the quasi-regularity,
the Dirichlet form is associated with a diffusion process on the Wasserstein space, which will be addressed in the forthcoming paper \cite{RWW}.

Let $\scr P_2$ be the space of all probability measures on $\R^d$ having finite second moment. It is a complete and separable   space under the quadratic Wasserstein distance
 $$\W_2(\mu,\nu):= \inf_{\pi\in\C(\mu,\nu)} \bigg(\int_{\R^d\times \R^d} |x-y|^2\pi(\d x,\d y) \bigg)^{\ff 1 2},\ \ \mu,\nu\in \scr P_2,$$
 where $\C(\mu,\nu)$ is the set of all couplings of $\mu$ and $\nu$. This space has been equipped with a natural Riemannian structure and becomes an infinite-dimensional Riemannian manifold, which is
 non-standard as explained in Remark 2.1(1) below. This structure is induced by the intrinsic derivative introduced by  Albeverio, Kondratiev and R\"ockner \cite{AKR}, and is consistent with
 Otto's structure   \cite{Otto} defined for probability measures having smooth and positive density functions,   see Remark  2.1 below for details. The Wasserstein space is a fundamental research object in the theory of optimal transport and related analysis, see \cite{AMB, V}
 and reference therein.  In this paper we aim to study stochastic analysis on the Wasserstein space. To this end, we first recall some existing literature concerning measure-valued diffusion processes.

Firstly,  to construct   measure-valued diffusion processes using the  Dirichlet-form  theory,   the integration by parts formula has been  established  for the intrinsic and extrinsic derivatives with respect to a reference distribution $\Xi$ on the space of Radon measures, see   \cite{1,KLV, ORS, C, Shao, Sturm} and references therein. Moreover, functional inequalities have been derived for measure-valued processes, see  \cite{DS, FMW, FSW, FW, FW2, RW19, ST, W18,  WZ}.  A key point in the construction of a Dirichlet form is to establish the integration by parts formula of the reference measure for derivatives in measure.
  The stationary distributions   in these references are chosen as either the  entropic-type measures supported on the space of singular distributions without discrete part,
 or the Dirichlet/Gamma type measures   concentrated on  the space of discrete distributions, which have reasonable backgrounds from physics,   population genetics and Bayesian non-parametrics.   Along a different direction,  a Rademacher type theorem is established in \cite{B}  for a class of reference measures  satisfying the integration by parts formula for the intrinsic derivative.

Next, corresponding to Dean-Kawasaki type SPDEs, local Dirichlet forms have been constructed on the   Wasserstein  space over the real line induced by increasing functions
(see \cite{Ko1, Ko2} and references therein), and it is proved in \cite{Ko3} that  the  associated diffusion process is given by
the empirical measure of independent particle systems.

Moreover, by solving  a conditional distribution-dependent SDE, \cite{W21} constructed a  diffusion process on $\scr P_2$ with generator  given by a second-order differential operator
in intrinsic derivative, and establish the Feynman$-$Kac formula for the underline measure-valued PDE. Since the SDE is driven by finite-dimensional Brownian motion,   the measure-valued diffusion process  is highly degenerate. Note  that the measure-valued diffusion processes constructed in \cite{W21}  extends that generated by the partial Laplacian investigated in \cite{A}.
 See also \cite{DFL} for an extension to the Wasserstein space over a compact Riemannian manifold.

 \

In this paper, we   construct and study  the Ornstein$-$Uhlenbeck (OU) type Dirichlet form  on the Wasserstein space $\scr P_2$, whose stationary distribution is  a fully supported Gaussian measure,     and the
generator is the Laplacian with a drift, where the Laplacian is induced by the Riemannian structure and is hence   crucial in geometric analysis.    In view of the fundamental role played by   Ornstein$-$Uhlenbeck process  in Malliavin calculus  on the Wiener  space, the present study should be crucial for developing stochastic analysis on the Wasserstein space.

Recall  that the Brownian motion on a $d$-dimensional Riemannian manifold can be constructed by using  the flat Brownian motion   on the tangent space $\R^d$. In the same spirit, we will   introduce the tangent space on $\scr P_2$ which is a separable Hilbert space, then recall the Gaussian measure and Ornstein$-$Uhlenbeck process on the Hilbert space, and finally construct the corresponding objects on $\scr P_2$ as projections from the tangent space.  The Ornstein$-$Uhlenbeck type Dirichlet form  we construct  on $\scr P_2$ shares nice properties
of the original process on the tangent space: it satisfies the log-Sobolev inequality and the generator has purely discrete spectrum.

  \

 The remainder of the paper is organized as follows.
In Section 2,  we  recall   the Riemannian structure induced by the intrinsic derivative due to \cite{AKR}, and calculate  the Laplacian operator $ \DD_{\scr P_2}.  $ Since different structures and Laplacians exist in the literature,
to distinguish we call ours intrinsic Riemannian structure and intrinsic Laplacian.
In Section 3,  we construct the Gaussian measure $N_{\mu_0,Q}$ determined  by a reference measure $\mu_0\in \scr P_2$ together with
 an unbounded positive definite linear operator $Q$ on the tangent space $T_{\mu_0}:=L^2(\R^d\to\R^d,\mu_0)$, and study the corresponding OU type Dirichlet form on $\scr P_2$.
   In Section 4, we try to formulate the generator of the OU process as
    $$Lf(\mu)= \DD_{\scr P_2} f(\mu)- \<b(\mu),  Df(\mu)\>_{T_\mu},$$
 where  $T_\mu:=L^2(\R^d\to\R^d,\mu)$ is the tangent space at $\mu$, and    $b(\mu)\in T_\mu$ is induced by the linear operator $Q$ on $T_{\mu_0}$. This
 formulation is consistent with that of the  OU process  on a separable Hilbert space, but it is in the moment $``$either formal computation or non-rigorous formula"
 due to the lack of a reasonable class of functions making the right hand side  meaningful.
  Finally, in Section 5 we study symmetric diffusion processes on $\scr P_2$ as perturbations of the OU process.

 \section{Intrinsic Riemannian structure on the Wasserstein space}

By using the gradient flow of density functions  arising from Monge's optimal transport, Otto \cite{Otto} constructed the Riemannian structure   on   $\scr P_2^{ac},$
   the  space of  measures in $\scr P_2$  having strictly positive smooth density functions with respect to the Lebesgue measure, see also  \cite[Chapter 13]{V}.
Under Otto's  structure,  the tangent space at $\mu\in \scr P_2^{ac}$ is the $L^2(\mu)$-closure of $\{\nn f: f\in C_b^\infty(\R^d)\}$. The Ricci curvature  was calculated in \cite{LOTT}, while the Levi-Civita connection and parallel displacement have been studied in \cite{DF}.

To built up a Riemannian structure, one needs to introduce the tangent space at each point $\mu\in \scr P_2$, and define an inner product (or metric) on the tangent space. From different point of views, several    Riemannian structures on the Wasserstein space  have been
studied in the literature, see \cite{E} for
 the tangent space  induced by all germs of geodesic curves, and see \cite{Otto} for
 the tangent space   induced by germs of geodesic curves induced by optimal transport maps (the closure of gradient-type vector fields).
In this paper, we adopt the $L^2$-tangent space   introduced in  \cite{AKR} (see \cite[Appendix]{RRW}) to define  the intrinsic derivative. This structure fits well to  the G\^ateaux derivative in infinite-dimensional analysis, and it works for the space of general Radon measures as well.


 \subsection{Intrinsic derivative}

We will simply denote $\mu(f)=\int f\d\mu$ for a measure $\mu$ and a function $f\in L^1(\mu)$.
For any $\mu\in \scr P_2$ and measurable $\phi:\R^d\to\R^d$, let
$\mu\circ\phi^{-1}$ be the image of $\mu$ under $\phi$, i.e.
$$(\mu\circ \phi^{-1} )(A):= \mu(\phi^{-1}(A))$$ for measurable sets $A\subset \R^d$. It is easy to see that $\mu\circ \phi^{-1}\in \scr P_2$ if and only if
$$\phi\in T_\mu:=L^2(\R^d\to\R^d,\mu),$$ where $L^2(\R^d\to\R^d,\mu)$ is the space of measurable maps $\phi$ from $\R^d$ to $\R^d$ with
$$\|\phi\|_{L^2(\mu)}:=(\mu(|\phi|^2))^{\ff 1 2}<\infty.$$ So, it is natural to take $T_\mu$ as the tangent space
at $\mu$, which is a separable Hilbert space with inner product
\beq\label{MT} \<\phi_1,\phi_2\>_{T_\mu}:= \mu\big(\<\phi_1,\phi_2\>\big)=\int_{\R^d} \<\phi_1,\phi_2\>\d\mu,\ \ \phi_1,\phi_2\in T_\mu.\end{equation}
Let $id\in T_\mu$ be the identity map, i.e. $id(x)=x$.

\beg{defn} Let $f\in C(\scr P_2)$, the class of continuous functions on $\scr P_2$.
\beg{enumerate} \item[$(1)$] We call $f$ intrinsically differentiable, if for any $\mu\in \scr P_2$,
$$T_\mu\ni \phi\mapsto D_\phi f(\mu):=\lim_{\vv\downarrow 0} \ff{f(\mu\circ (id+\vv \phi)^{-1})- f(\mu)}\vv\in\R $$
is a   bounded linear functional. In this case, the intrinsic derivative $Df(\mu)$ is the unique element in $T_\mu$ such that
$$\<Df(\mu),\phi\>_{T_\mu}:=\mu(\<\phi, Df(\mu)\>)= D_\phi f(\mu),\ \ \phi\in T_\mu.$$
\item[$(2)$]   We write $f\in  C^1(\scr P_2)$, if $f$ is  intrinsically differentiable  such that $Df(\mu)(x)$ has a  jointly continuous version;
 i.e. each $Df(\mu)$ has a $\mu$-version $x\mapsto Df(\mu)(x)$ such that
 $$(\mu,x)\in \scr P_2\times \R^d\mapsto Df(\mu)(x)$$ is continuous. We denote $f\in C_b^1(\scr P_2)$ if moreover   $f$ and $Df$ are bounded.
\item[$(3)$] We denote $f\in  C^2(\scr P_2)$, if $f\in C^1(\scr P_2)$,  the jointly continuous version  $Df(\mu)(x)$ is intrinsically differentiable in $\mu$ and differentiable in $x$, such that
$$D^2f(\mu)(x,y):=  D\big(Df(\cdot)(x)\big)(\mu)(y),\ \ \nn Df(\mu)(x):= \nn\big(Df(\mu)(\cdot)\big)(x)$$ have versions     jointly continuous in all arguments in the same sense of item (2).  We write $f\in C_b^2(\scr P_2)$ if moreover $f, Df, D^2f$ and $\nn Df$ are bounded. \end{enumerate}
\end{defn}

When $f\in C^1(\scr P)$, we automatically take $Df(\mu)(x)$ to be the jointly continuous version of $Df$, which is unique.  Indeed, by the continuity,
$Df(\mu)(\cdot)$ is unique for each $\mu\in \scr P_2$ with full support, so that it is unique for all $\mu\in \scr P_2$ since the set of fully supported  measures is dense in $\scr P_2$.
Under the Riemannian metric  given by  \eqref{MT},
 the space $\scr P_2$ becomes an infinite-dimensional $``$Riemannian manifold with different level  boundaries," see Remark 2.1(1) below for an explanation.    

To make calculus on $\scr P_2$,  we   introduce the displacement of the tangent space.
For any $\phi\in T_\mu$, consider the displacement of measures   along $\phi$ from $\mu$:
$$[0,\infty)\ni s\mapsto \mu\circ(id+s\phi)^{-1}\in \scr P_2.$$ Then the  tangent space is shifted as
\beq\label{TG} T_{\mu\circ (id+s \phi)^{-1}}= T_\mu \circ(id+s\phi)^{-1}:= \big\{h\circ(id+s\phi)^{-1}:\ h\in T_\mu\big\},\ \ s\ge 0,\end{equation}
where $h\circ(id+s\phi)^{-1}\in T_{\mu \circ(id+s\phi)^{-1}} $ is uniquely determined by
\beq\label{PO}\big\<h\circ(id+s\phi)^{-1}, \psi \big\>_{T_{ \mu \circ(id+s\phi)^{-1}}}:= \big\<h, \psi\circ(id+ s\phi)\big\>_{T_\mu},\ \ \psi\in T_{\mu \circ(id+s\phi)^{-1}},\end{equation}
by noting that   $\psi\circ(id+ s\phi)\in T_\mu$  is due to
\beq\label{PPS} \beg{split}& \|\psi\circ(id+ s\phi)\|_{T_\mu}^2= \mu(|\psi\circ(id+s\phi)|^2)\\
&=\big(\mu\circ(id+s\phi)^{-1}\big)(|\psi|^2) <\infty,\ \ \ \ \psi\in T_{\mu\circ(id+s\phi)^{-1}}.\end{split} \end{equation}
Obviously,    $T_{\mu\circ (id+s \phi)^{-1}}\supset T_\mu \circ(id+s\phi)^{-1}.$
 On the other hand,
for any $\psi\in T_{\mu\circ (id+s \phi)^{-1}}$, \eqref{PPS} implies
  $\tt\psi:=\psi\circ(id+s\phi)\in T_\mu$ and
$$\psi= \tt\psi\circ(id+s\phi)^{-1} \in T_\mu\circ (id+s\phi)^{-1}.$$
Therefore, \eqref{TG} holds.

The following result implies  that a function $f\in C_b^1(\scr P_2)$  is $L$-differentiable, i.e. it is intrinsically differentiable and
\beq\label{L} \lim_{\|\phi\|_{L^2(\mu)}\downarrow 0} \ff{|f(\mu\circ(id+\phi)^{-1})-f(\mu)-D_\phi f(\mu)|}{\|\phi\|_{L^2(\mu)}} =0.\end{equation}
In this case, the intrinsic derivative is also called the $L$-derivative, which   coincides with Lions' derivative introduced in   \cite{Card}.

\beg{prp}\label{P} Let $f\in C^1(\scr P_2)$ such that
\beq\label{UI}\lim_{N\to\infty} \limsup_{\|\phi\|_{L^2(\mu)}\downarrow 0}       \Big\|\big(|Df(\mu\circ(id+\phi)^{-1})(id+\phi)|-N\big)^+\Big\|_{L^2(\mu)}
=0,  \end{equation}  then $\eqref{L}$ holds, i.e. $f$ is $L$-differentiable. \end{prp}

\begin{proof} Let $\mu\in \scr P_2$ and $\phi\in T_\mu$. By \eqref{TG} we have $\phi\circ (id+s\phi)^{-1} \in T_{\mu\circ(id+s\phi)^{-1}}$ for $s\in [0,1]$, and
\beg{align*} &\ff{\d}{\d s} f(\mu\circ (id+s\phi)^{-1})=\lim_{\vv\downarrow 0} \ff{f(\mu\circ (id+(s+\vv)\phi)^{-1}) - f(\mu\circ(id+s\phi)^{-1})}\vv \\
&= \lim_{\vv\downarrow 0} \ff{f\big((\mu\circ (id+s \phi)^{-1}) \circ(id+\vv \phi\circ(id+s\phi)^{-1})^{-1} \big) - f(\mu\circ(id+s\phi)^{-1})}\vv\\
&= D_{\phi\circ(id+s\phi)^{-1}} f(\mu\circ(id+s\phi)^{-1})= \mu\big(\<\phi, (Df(\mu\circ(id+s\phi)^{-1}))\circ(id+s\phi)\>\big).\end{align*}
Combining this with  $f\in C^1(\scr P_2)$,   we arrive at
\beg{align*} &\limsup_{\|\phi\|_{L^2(\mu)}\downarrow 0} \ff{|f(\mu\circ(id+\phi)^{-1}) - f(\mu)-D_\phi f(\mu)|}{\|\phi\|_{L^2(\mu)}}\\
&=\limsup_{\|\phi\|_{L^2(\mu)}\downarrow 0} \ff{|\int_0^1 \ff{\d}{\d s} f(\mu\circ (id+s\phi)^{-1}) \d s -D_\phi f(\mu)|}{\|\phi\|_{L^2(\mu)}}\\
&\le  \limsup_{\|\phi\|_{L^2(\mu)}\downarrow 0} \int_0^1 \ff{ 1}{\|\phi\|_{L^2(\mu)}} \Big|\mu\big(\<\phi, (Df(\mu\circ(id+s\phi)^{-1}))\circ(id+s\phi)- Df(\mu)\>\big)\Big| \d s\\
&\le  \limsup_{\|\phi\|_{L^2(\mu)}\downarrow 0} \int_0^1 \big\| (Df(\mu\circ(id+s\phi)^{-1}))\circ(id+s\phi)- Df(\mu)\big\|_{L^2(\mu)} \d s =0,\end{align*}
where the last step follows from the continuity of $Df$,    \eqref{UI} and the dominated convergence theorem.
Indeed, if the last step does not hold, then there exist a constant $\vv>0$ and a sequence $\{\phi_n\}_{n\ge 1} $ with $\|\phi_n\|_{L^2(\mu)}\le \ff 1 n$ such that
$$\xi_n(s):= Df(\mu\circ(id+s\phi_n)^{-1})(id+s\phi_n),\ \ n\ge 1, s\in [0,1]$$ satisfies
\beq\label{71} \int_0^1\|\xi_n(s)- Df(\mu)\|_{L^2(\mu)}\d s\ge \vv,\ \ n\ge 1.\end{equation}
Up to a subsequence, we may also assume that $\mu(\{\phi_n\to 0\})=1$, which together with
$$\sup_{s\in [0,1]}\W_2(\mu\circ(id+s\phi_n)^{-1}, \mu)\le \|\phi_n\|_{L^2(\mu)}\le \ff 1 n,\ \ n\ge 1 $$ and $f\in C^1(\scr P_2)$    implies
$$\mu\Big(\lim_{n\to\infty} \xi_n(s)= Df(\mu), \ s\in [0,1]\Big)=0.$$
So, by   the dominated convergence theorem,  we obtain
\beq\label{72} \lim_{n\to\infty} \int_0^1 \big\|(\xi_{n}(s)-Df(\mu)) 1_{\{|\xi_n(s)|<N\}}\big\|_{L^2(\mu)} \d s=0,\ \ N\ge 1.\end{equation}
On the other hand, \eqref{UI} implies
\beg{align*} &\lim_{N\to\infty} \limsup_{n\to\infty}  \int_0^1\big\|(|\xi_n(s)|+|Df(\mu)|) 1_{\{|\xi_n(s)|\ge 2 N\}}\big\|_{L^2(\mu)}  \d s\\
&\le \lim_{N\to\infty} \limsup_{n\to\infty}  \int_0^1\big\|2 (|\xi_n(s)|+|Df(\mu)|  -N )^+   \big\|_{L^2(\mu)}  \d s\\
&\le \lim_{N\to\infty} \limsup_{\|\phi\|_{L^2(\mu)}\downarrow 0}   2  \Big\|\big(|Df(\mu\circ(id+\phi)^{-1})(id+\phi)|+|Df(\mu)| -N \big)^+\Big\|_{L^2(\mu)}  =0. \end{align*}
Combining this with \eqref{72} leads to
\beg{align*} &\limsup_{n\to\infty}  \int_0^1\|\xi_n(s)- Df(\mu)\|_{L^2(\mu)}\d s\\
&\le \limsup_{N\to\infty}\limsup_{n\to\infty}     \int_0^1 \big\|(\xi_n(s)- Df(\mu)) 1_{\{|\xi_n(s)|< 2 N\}}\big\|_{L^2(\mu)}\d s\\
&\quad
 +\limsup_{N\to\infty}\limsup_{n\to\infty}      \int_0^1 \big\|(|\xi_n(s)|+|Df(\mu)|) 1_{\{|\xi_n(s)|\ge 2 N\}}\big\|_{L^2(\mu)}   \d s  =0,\end{align*}
which contradicts to \eqref{71}.

\end{proof}

We are ready to introduce the chain rule  for the intrinsic derivative  in the distribution   of random variables. Let $\L_\xi$ be the law of a random variable under a probability space
$(\OO,\F,\P)$.  Note that \eqref{09} implies
\beg{align*}&\big(|Df(\mu\circ(id+\phi)^{-1})(id+\phi)|-N\big)^+\\
&\le \big(c(1+|x|+ |\phi(x)|)-N\big)^+\le c|\phi(x)|+ \big(c +c|x| -N\big)^+,\end{align*}
so that \eqref{UI} holds for any $\mu\in \scr P_2.$
Then
the following result follows from Proposition \ref{P} and \cite[Theorem 2.1(2)]{BRW} for $p=2$, see also \cite[Lemma A.8]{18} for an earlier result.

\beg{prp}\label{P1} Let $(\xi_s)_{s\in [0,1]}$ be a family of $\R^d$-valued random variables on a probability space $(\OO,\F,\P)$,  such that $\mu_s:=\L_{\xi_s}\in \scr P_2$ and
$$\dot \xi_0:=\lim_{s\downarrow 0} \ff{\xi_s-\xi_0}s$$ exists in $L^2(\P)$.  Then for any $f\in C^1(\scr P_2)$ such that
\beq\label{09} |D f(\mu)(x)|\le c(1+|x|),\ \ x\in\R^d,\ \ \mu\in \scr P_2, \W_2(\mu,\mu_0)\le 1\end{equation}  holds for some constant $c>0$, we have
$$\ff{\d}{\d s}\Big|_{s=0} f(\L_{\xi_s}):=\lim_{s\downarrow 0} \ff{f(\L_{\xi_s})- f (\L_{\xi_0})}s= \E\big[\<Df(\mu_0)(\xi_0), \dot\xi_0\>\big].$$
\end{prp}
\ \newline
{\bf Remark 2.1.}   Below we make some comments on  the above Riemannian structure.

\beg{enumerate} \item[(1)] Since the dimension of the tangent space $T_\mu$ is $d$ times the number  of points contained in the support of $\mu$, tangent spaces at different points may have different dimensions. This is different from
a standard Riemannian manifold, but fits well to the feature of Riemannian manifolds with different level boundaries (sub-manifolds).
For instance,  consider a half ball in $\R^3$. It is a 3D manifold  with 2D boundary, and the 2D boundary itself is a manifold with 1D boundary.
At a point on the boundary outside the edge, the tangent space only consists of vectors tangent to the boundary, since the geodesic along the normal vector may go beyond the manifold, so geodesics   from this point are only defined along tangent vectors of the boundary  under the induced Riemannian metric, so that the Laplacian becomes  degenerate comparing to 3D, and the Brownian motion starting from this point stays on the 2D boundary. The same happens to points on the 1D boundary.

 \item[(2)]  The intrinsic derivative coincides with Otto's derivative when  it exists, where the latter is only  defined for absolutely continuous measures in $\scr P_2$.
See  \cite{E} for an extension of Otto's structure to the whole space  $\scr P_2$.
\item[$(3)$] To see that $\W_2$ is the intrinsic distance of the Riemannian structure,   let $\mu_0\in \scr P_2$ be absolutely continuous with respect to the Lebesgue measure.
Then for any   $\mu_1,\mu_2\in \scr P_2,$  there exists
an optimal coupling
$(h_1,h_2)\in T_{\mu_0}\times T_{\mu_0},$ in the sense of random variables on $\R^d$ under the probability measure $\mu_0$,  such that
$$\mu_i=\mu_0\circ h_i^{-1} (i=1,2),\ \ \W_2(\mu_1,\mu_2)^2=\mu_0(|h_1-h_2|^2),$$
so that
\beq\label{NT} \nu_t:= \mu_0\circ(th_1+(1-t)h_2)^{-1},\ \ \ t\in [0,1]\end{equation}  is the geodesic linking $\mu_1$ and $\mu_2$, i.e.
$$\nu_0=\mu_2, \ \ \nu_1=\mu_1,\ \ \W_2(\nu_s,\nu_t)= |t-s|\W_2(\mu_1,\mu_2)\ \text{for}\ t,s\in [0,1].$$
Indeed, for any $0\le s\le t\le 1$, \eqref{NT} implies
$$\mu_0\circ \psi_{s,t}^{-1}\in \C(\nu_s,\nu_t)$$ for $$  \psi_{s,t}:= \big(sh_1+(1-s)h_2,  th_1+(1-t)h_2\big): \R^d\to \R^d\times\R^d,$$
so that by the definition  of $\W_2$, we have
\beg{align*} &\W_2(\nu_s,\nu_t)\le \bigg(\int_{\R^d\times \R^d} |x-y|^2 (\mu_0\circ \psi_{s,t}^{-1})(\d x,d y)\bigg)^{\ff 1 2}\\
& = \bigg(\int_{\R^d }  |th_1+(1-t)h_2-sh_1-(1-s)h_2|^2(x)\mu_0(\d x)\bigg)^{\ff 1 2} \\
&\le |t-s| \ss{\mu_0(|h_1-h_2|^2)}=(t-s)\W_2(\mu_1,\mu_2),\ \ \ \ 0\le s\le t\le 1,\end{align*}
which together with   the triangle inequality implies
$$\W_2(\mu_1,\mu_2)\le \W_2(\nu_0, \nu_s)+\W_2(\nu_s,\nu_t)+\W_2(\nu_t,\nu_1) \le \W_2(\mu_1,\mu_2),\ \ 0\le s\le t\le 1,$$
and hence,  $ \W_2(\nu_s,\nu_t)= |t-s|\W_2(\mu_1,\mu_2)$   for all $t,s\in [0,1]$.\end{enumerate}

\subsection{Intrinsic Laplacian}

Recall that on a $d$-dimensional Riemannian manifold, the Laplacian is defined as the trace of the second-order derivative  $\nn^2$ (i.e. Hessian operator). Below we define the intrinsic Laplacian on $\scr P_2$ in the same way.

Let $\mu\in \scr P_2$ and let $\{\phi_m\}_{m\ge 1}$ be an ONB (orthonormal basis) of $T_\mu:= L^2(\R^d\to\R^d,\mu)$. Note that  the number of $\{\phi_m\}_{m\ge 1}$ is a finite family if and only if
$\mu$ has finite support. For any $\phi\in T_\mu,$ let
$$D_\phi^2 f(\mu):=
 \ff{\d}{\d\vv}\Big|_{\vv=0} \big(D_{\phi\circ (id+\vv \phi)^{-1}} f\big)(\mu\circ(id+\vv \phi)^{-1}),\ \ f\in C_b^2(\scr P_2).$$
We write
$f\in \D_\mu(\DD_{\scr P_2}),$  the domain of $\DD_{\scr P_2}$ at point $\mu$, if $f\in C_b^2(\scr P_2)$ is such that
\beq\label{LA1} \beg{split} &\DD_{\scr P_2} f(\mu):= \sum_{m\ge 1} D_{\phi_m}^2 f(\mu)  \end{split}\end{equation}
exists.   Note that   the partial Laplacian considered in \cite{A} is given by
$$\DD_w f(\mu):= \sum_{i=1}^dD^2_{e_i} f(\mu),$$
where $\{e_i\}_{1\le i\le d}$ is the standard orthonormal basis of $\R^d$, which is thus a partial sum from that in the definition of $\DD_{\scr P_2} f(\mu).$ Since $\mu$ is a probability measure, the constant vectors $\{e_i\}_{1\le i\le d}$ are orthonormal in $T_\mu\scr P_2$.
We have the following formulation of $\DD_{\scr P_2}$.

\beg{prp}\label{LP} For any $\mu\in \scr P_2$ and $f\in \D_\mu(\DD_{\scr P_2})$,
  \beg{align*} \DD_{\scr P_2} f(\mu)=  \sum_{m\ge 1}  \int_{\R^d\times\R^d} \Big(&\<D^2 f(\mu)(x,y), \phi_m(x)\otimes \phi_m(y)\>_{HS}\\
  & +\<\nn D f(\mu)(x), \phi_m(x)\otimes \phi_m(x)\>_{HS}\Big)\mu(\d x)\mu(\d y), \end{align*}
 where $\<\cdot,\cdot\>_{HS}$ is the Hilbert$-$Schmidt inner product for matrices, and the right-hand side   does not depend on the choice of the ONB $\{\phi_m\}_{m\ge 1}.$
 Consequently,  for any $f\in C_b^2(\scr P_2)$ and $\mu\in \scr P_2$, we have  $f \in \D_\mu(\DD_{\scr P_2})$ if and only if the following series exists:
$$ {\rm tr}\big(\nn Df(\mu)\big):= \sum_{m\ge 1}\int_{\R^d} \<\nn D f(\mu)(x), \phi_m(x)\otimes \phi_m(x)\>_{HS}\, \mu(\d x).$$
 \end{prp}
\beg{proof}
Noting that
\beg{align*} &\big(D_{\phi_m\circ (id+\vv \phi_m)^{-1} }f\big)(\mu\circ(id+\vv \phi_m)^{-1})\\
&=\big(\mu\circ(id+\vv\phi_m)^{-1} \big)\big(\<\phi_m\circ(id+\vv\phi_m)^{-1}, Df(\mu\circ (id+\vv\phi_m)^{-1})\>\big)\\
&=\mu\big(\<\phi_m, Df(\mu\circ (id+\vv \phi_m)^{-1})(id+\vv\phi_m)\>\big),\end{align*}
we obtain
\beg{align*} &D_{\phi_m}^2 f(\mu)=\ff{\d}{\d\vv}\Big|_{\vv=0} \mu\big(\<\phi_m, Df(\mu\circ (id+\vv \phi_m)^{-1})(id+\vv\phi_m)\>\big)  \\
&= \int_{\R^d\times\R^d} \<D^2 f(\mu)(x,y), \phi_m(x)\otimes \phi_m(y)\>_{HS}\,\mu(\d x)\mu(\d y)\\
&\qquad + \int_{\R^d} \<\nn D f(\mu)(x), \phi_m(x)\otimes \phi_m(x)\>_{HS}\, \mu(\d x),\end{align*}
where the first term comes from the derivative  $\ff{\d}{\d\vv}|_{\vv=0} Df(\mu\circ (id+\vv \phi_m)^{-1})$ by Proposition \ref{P1} for $\xi_\vv:=  id+\vv \phi_m,$ and the other term
follows from the derivative  $ \ff{\d }{\d\vv}|_{\vv=0} Df(\mu)(id+\vv\phi_m).$
 Therefore, the desired formula holds by \eqref{LA1}.

 Next, it is easy to see that $\DD_{\scr P_2} f(\mu)$  does not depend on the choice of the ONB
$\{\phi_m\}_{m\ge 1}$.   Indeed, for another ONB  $\{\tt \phi_m\}_{m\ge 1}$   in $T_\mu$, we have
\beg{align*} &\tt \phi_m= \sum_{l\ge 1} \<\tt\phi_m,\phi_l\>_{T_\mu} \phi_l,\ \ \ m\ge 1,\\
&\sum_{m\ge 1} \<\tt\phi_m, \phi_l\>_{T_\mu}\<\tt\phi_m, \phi_k\>_{T_\mu}= 1_{\{l=k\}},\ \ k,l\ge 1.\end{align*}
So,
\beg{align*} &\sum_{m\ge 1}\Big(\<D^2 f(\mu)(x,y), \tt\phi_m(x)\otimes \tt\phi_m(y)\>_{HS}+\<\nn D f(\mu)(x), \tt\phi_m(x)\otimes \tt\phi_m(x)\>_{HS}\Big)\\
&=\sum_{k,l,m\ge 1} \<\tt\phi_m, \phi_k\>_{T_\mu} \<\tt\phi_m, \phi_l\>_{T_\mu}\\
&\qquad\qquad \times  \Big(\<D^2 f(\mu)(x,y), \phi_k(x)\otimes \phi_l(y)\>_{HS}
 +\<\nn D f(\mu)(x), \phi_k(x)\otimes \phi_l(x)\>_{HS}\Big)\\
&= \sum_{k,l\ge 1} 1_{\{k=l\}}  \Big(\<D^2 f(\mu)(x,y), \phi_k(x)\otimes \phi_l(y)\>_{HS}+\<\nn D f(\mu)(x), \phi_k(x)\otimes \phi_l(x)\>_{HS}\Big)\\
&=\sum_{m=1}\Big(\<D^2 f(\mu)(x,y), \phi_m(x)\otimes \phi_m(y)\>_{HS}+\<\nn D f(\mu)(x), \phi_m(x)\otimes \phi_m(x)\>_{HS}\Big).\end{align*}

Finally, let   $\{e_i\}_{1\le i\le d}$ be the standard ONB in $\R^d$. For any $f\in C_b^2(\scr P_2)$, by the Cauchy$-$Schwarz inequality  we obtain
\beg{align*} &\sum_{m\ge 1} \bigg|\int_{\R^d\times \R^d} \<D^2 f(\mu)(x,y), \phi_m(x)\otimes \phi_m(y)\>_{HS}\,\mu(\d x)\mu(\d y) \bigg|\\
&= \sum_{m\ge 1} \bigg|\sum_{i,j=1}^d \int_{\R^d\times \R^d} \big(D^2 f(\mu)(x,y)\big)_{ij}  \<e_i, \phi_m(x)\>\<e_j,   \phi_m(y)\>\,\mu(\d x)\mu(\d y) \bigg|\\
&\le \|D^2f\|_\infty \sum_{i,j=1}^d \bigg(\sum_{m\ge 1} \<e_i, \phi_m\>_{T_\mu}^2\bigg)^{\ff 1 2} \bigg(\sum_{m\ge 1} \<e_j, \phi_m\>_{T_\mu}^2\bigg)^{\ff 1 2} \\
&= \|D^2f\|_\infty \sum_{i,j=1}^d \|e_i\|_{L^2(\mu)}\|e_j\|_{L^2(\mu)}= d^2 \|D^2 f(\mu)\|_\infty<\infty.\end{align*}
Therefore, $ f\in \D_\mu(\DD_{\scr P_2})$ if and only if  $ {\rm tr}\big(\nn Df(\mu)\big)$ exists.

 \end{proof}
\ \newline
{\bf Remark 2.2.} We present some comments on $\D_\mu(\DD_{\scr P_2})$. \beg{enumerate}
\item[$(1)$]  If $ \mu$ has finite support, then $T_\mu$ is finite-dimensional so that $ C_b^2(\scr P_2)\subset \D_\mu(\DD_{\scr P_2}).$
\item[$(2)$] Let  $f_i(\mu):=\mu(\<\cdot,e_i\>)= \int_{\R^d}x_i\mu(\d x), 1\le i\le d$, where    $\{e_i\}_{1\le i\le d}$ is  the standard ONB in $\R^d$ so that
$\<x,e_i\>=x_i.$  Then for any $\mu\in \scr P_2$ and $g\in C_b^2(\R^d)$,
$$f:=g(f_1,\dotsc, f_d)\in \D_\mu(\DD_{\scr P_2}).$$  Indeed,  this case we have
\beg{align*} &Df(\mu)=(\nn   g)(f_1(\mu),\dotsc, f_d(\mu)),\\
&D^2 f(\mu)= (\nn^2 g) (f_1(\mu),\dotsc, f_d(\mu)),\end{align*}
so that $\nn Df=0$ and for the standard ONB   $\{e_i\}_{1\le i\le d}$ in  $\R^d$,
\beg{align*} &\DD_{\scr P_2} f(\mu)=\sum_{m\ge 1} \sum_{i,j=1}^d (\pp_i\pp_j g)(f_1(\mu),\dotsc, f_d(\mu)) \int_{\R^d\times\R^d} \phi_m^i(x)\phi_m^j(y)\mu(\d x)\mu(\d y) \\
&=   \sum_{i,j=1}^d (\pp_i\pp_j g)(f_1(\mu),\dotsc, f_d(\mu))  \sum_{m=1}^d \mu(\<\phi_m, e_i\>)\mu(\<\phi_m, e_j\>) \\
&= \sum_{i,j=1}^d (\pp_i\pp_j g)(f_1(\mu),\dotsc, f_d(\mu))\mu(\<e_i, e_j\>)\\
&= (\DD g)(f_1(\mu),\dotsc, f_d(\mu)),\ \ \mu\in \scr P_2.\end{align*}
\item[$(3)$]  In general,  non-constant cylindrical functions of type
$$f(\mu):= g(\mu(f_1),\dotsc, \mu(f_n)),\ \ n\ge 1, f_i\in C_b^2(\R^d), g\in C_b^2(\R^n)$$
are not in  $\D(\DD_{\scr P_2}):=\cap_{\mu\in \scr P_2} \D_\mu(\DD_{\scr P_2}),$ although they are in $C_b^2(\scr P_2),$ see    \cite[Remark  5.6]{E} for explanations.
 So,    the domain of $\DD_{\scr P_2}$ might be too small to define a Brownian motion type diffusion process on $\scr P_2$. In   the next section,  we
 construct   OU type Dirichlet forms  on $\scr P_2$, so that   the associated  generators have dense domains in the $L^2$-space,   although the structure of functions in the domains is  unknown, see Remark 4.1 below.

\end{enumerate}

 \section{Ornstein$-$Uhlenbeck process on Wasserstein space}

We will start from the OU process on  the tangent space $T_{\mu_0}$ for a reference measure $\mu_0\in \scr P_2$, then transform to the Wasserstein space.

 To make sure that
any $\mu\in \scr P_2$ is the distribution of some $h\in T_{\mu_0}$ under the probability $\mu_0$, i.e. $\mu=  \mu_0\circ h^{-1}$, we assume that $\mu_0$ is absolutely continuous
with respect to the Lebesgue measure.
 In this case,  for any $\mu\in \scr P_2$, there exists a unique $h\in T_{\mu_0}$ such that
$$\Psi(h):= \mu_0\circ h^{-1}=\mu,\ \ \ \ \W_2(\mu_0,\mu)^2= \mu_0(|id-h|^2).$$  This $h$ is called the optimal map as solution of the Monge problem for $\W_2$, see \cite[Theorem 10.41]{V} or \cite{AMB}.
 The map   $\Psi: T_{\mu_0}\to \scr P_2$ is a Lipschitz surjection, i.e. $\Psi  (T_{\mu_0} ) =\scr P_2$ and
 $$\W_2(\Psi(h),\Psi(\tt h))\le \mu_0(|h-\tt h|^2)^{\ff 1 2}= \|h-\tt h\|_{T_{\mu_0}},\ \ \ \ h,\tt h\in {T_{\mu_0}}.$$
 In the following, we first introduce some facts for the OU process on the Hilbert space $T_{\mu_0}$, then construct the corresponding one on $\scr P_2$.

\subsection{OU process on tangent   space}

Since the tangent space $T_{\mu_0}$ is a separable Hilbert space, it has a natural flat Riemannian structure. Let  $\{h_n\}_{n\ge 1}$  be a complete orthonormal basis of
$T_{\mu_0}$. The  Laplacian is given by
$$\DD := {\rm tr}(\nn^2)=\sum_{n=1}^\infty \nn_{h_n}^2,  $$
where $\nn_{h_n}$ is the directional derivative along $h_n$, and   the gradient operator $\nn$ is determined by
$$\<\nn f(h), h_n\>_{T_{\mu_0} }:=  \nn_{h_n}  f(h),\ \ n\ge 1, h\in T_{\mu_0}, f\in C^1(T_{\mu_0}).$$

  Let $Q$ be a positive definite unbounded self-adjoint operator in ${T_{\mu_0}}$ with discrete spectrum  $\{q_n\}_{n\ge 1}$ and eigenbasis $\{h_n\}_{n\ge 1}$ such that $q_n\uparrow\infty$ as $n\uparrow\infty$ and
 $$\sum_{n=1}^\infty q_n^{-1}<\infty.$$ So, $Q$ has trace-class inverse $Q^{-1}$ and
  the centred Gaussian measure on ${T_{\mu_0}}$ with covariance $Q^{-1}$ is given by (see \cite{Bog})
 $$G_Q(\d h):= \prod_{n=1}^\infty \Big(\ff{q_n}{2\pi}\Big)^{\ff 1 2} \exp\Big[-\ff{q_n\<h,h_n\>_{T_{\mu_0}}^2}{2}\Big]\d \<h,h_n\>_{T_{\mu_0}}$$
 under the coordinates $\{\<h,h_n\>_{T_\mu}\}_{n\ge 1}$ referring to the expansion
 $$h= \sum_{n=1}^\infty \<h,h_n\>_{T_{\mu_0}} h_n,\ \ \ h\in T_{\mu_0}.$$
 The associated OU process can be constructed as  (see \cite[(5.2.9) or (6.2.1)]{DP})
\beq\label{HT} h_t= \e^{-tQ } h_0 +\ss 2\int_0^t \e^{-(t-s)Q} \d W_s,\ \ t\ge 0,\end{equation}
where $W_t$ is the cylindrical Brownian motion on $T_{\mu_0}$, i.e.
 $$W_t=\sum_{n=1}^\infty B_t^n h_n,\ \ t\ge 0$$
 for independent one-dimensional Brownian motions $\{B_t^n\}_{n\ge 1}.$

Let  $(\tt L, \D(\tt L))$  be   generator of the OU process associated with $G_Q$, which is a negative definite self-adjoint operator in $L^2(G_Q)$ with domain  $\D(\tt L)$ including  the class of cylindrical functions
$\F C_b^2(T_{\mu_0})$ consisting of
$$h\mapsto \tt f(h):=F(\<h,h_1\>_{T_{\mu_0}},\dotsc, \<h,h_n\>_{T_{\mu_0}}),\ \ n\ge 1, F\in C_b^2(\R^n),$$
and satisfying, for any such  function,
\beg{align*}&\tt L \tt f(h) = \DD\tt f(h)-\<Q \nn f(h), h\>_{T_{\mu_0}}\\
&= \sum_{i=1}^n  \big(\pp_i^2   -q_i \<h,h_i\>_{T_{\mu_0}}\pp_i  \big) F (\<h,h_1\>_{T_{\mu_0}},\dotsc, \<h,h_n\>_{T_{\mu_0}}),\end{align*}
where $\DD$ and $\nn$ are the Laplacian and gradient operators on $T_{\mu_0}$ respectively. Moreover, the  integration by parts formula yields
\beq\label{DER} \tt \EE(\tt f,\tt g):=\int_{T_{\mu_0}} \<\nn\tt f,\nn\tt g\>_{T_{\mu_0}} \d G_Q=-\int_{T_{\mu_0}}( \tt f \tt L\tt g)\d G_Q,\ \ \ \tt f,\tt g\in \scr FC_b^2(T_{\mu_0}).\end{equation}
Consequently, $(\tt\EE,\scr FC_b^2(T_{\mu_0}))$ is closable in $L^2(G_Q)$ and the closure $(\tt\EE,\D(\tt\EE))$ is a symmetric conservative Dirichlet form, see for instance \cite{MR}.
Moreover, it  satisfies the log-Sobolev inequality (see \cite{G1,G2})
\beq\label{LS1} G_Q(\tt f^2\log\tt f^2)\le \ff 2 {q_1}  \tt\EE(\tt f,\tt f),\ \ \tt f\in \D(\tt\EE),\ G_Q(\tt f^2)=1. \end{equation}

The generator  $(\tt L, \D(\tt L))$ of $(\tt E,\D(\tt E))$ has  purely discrete spectrum,  i.e. its essential spectrum  is empty. Indeed, consider the following one-dimensional  OU operators
$\{L_i\}_{i\ge 1}$:
 $$L_i \varphi(r)=\varphi''(r)-q_i r \varphi'(r),\ \   \ \ r\in \R.$$
 It is well known that each $-L_i$ has purely discrete spectrum consisting of simple eigenvalues
 $$\si(-L_i)  =\{\ll_{i,k}: k\ge 0\},$$where  $\ll_{i,0}=0, \ll_{i,1}=q_i$ and $\ll_{i,k}\uparrow\infty$ with linear growth as $k\uparrow \infty$. Since $q_i\uparrow\infty$ as $i\uparrow\infty$ and $L$ is the sum of these operators, i.e.
 $$\tt Lf(h)=\sum_{i=1}^\infty L_i f_{i, h}(\<h,h_i\>_{T_{\mu_0}}),\ \ f_{i,h}(r):= \tt f(h-\<h,h_i\>_{T_{\mu_0}}h_i+ r h_i),$$
 the spectrum of $-\tt L$ is purely discrete with eigenvalues
 $$\sum_{i=1}^n \ll_{i, k_i},\ \ n\ge 1, k_i\ge 0.$$
 According to the spectral theory, the pure discreteness of the  spectrum for $\tt L$ is equivalent to the compactness  in $L^2(G_Q)$  of   the associated  Markov semigroup $\tt P_t:=\e^{\tt Lt}$
 for $t>0$, they are also equivalent to the compactness  in $L^2(G_Q)$ of the set
 $$ \big\{\tt f\in \D(\tt \EE):\ \tt\EE_1(\tt f):= \tt\EE(f,f)+G_Q(\tt f^2)\le 1\big\}.$$

Let  $C_b^1(T_{\mu_0})$ be  the class of all bounded  functions on ${T_{\mu_0}}$ with bounded and continuous Fr\'echet derivative.
By an approximation argument,  see  the proof of Lemma \ref{LMM} below for $F=0$, we have
   $\D(\tt\EE)\supset C_b^1(T_{\mu_0})$ and \eqref{DER} implies
 $$ \tt\EE(\tt f,\tt g)=  G_Q\big(\<\nn\tt f,\nn\tt g\>_{T_{\mu_0}}\big):=\int_{T_{\mu_0}} \<\nn\tt f,\nn\tt g\>_{T_{\mu_0}}\d G_Q,\ \ \tt f,\tt g\in C_b^1({T_{\mu_0}}).$$

\subsection{OU process on $\scr P_2$}
We first introduce the Gaussian measure and the corresponding OU Dirichlet form on $\scr P_2$.

\beg{defn} Let $\Psi:{T_{\mu_0}}\to\scr P_2, \Psi(h):=\mu_0\circ h^{-1}.$  \beg{enumerate} \item[$(1)$]
 $N_{\mu_0,Q}:= G_Q\circ\Psi^{-1}$
 is called the Gaussian measure on $\scr P_2$ with parameter $(\mu_0,Q)$.
\item[$(2)$]  Define the following OU bilinear form  on $L^2(N_{\mu_0,Q}):$
\beg{align*}
&  \D(\EE):= \big\{f\in L^2(N_{\mu_0,Q}):\ f\circ\Psi\in \D(\tt\EE)\big\},\\
&\EE(f,g):=\tt\EE(f\circ\Psi, g\circ\Psi),\  \ \ \ f,g\in\D(\EE).\end{align*}\end{enumerate} \end{defn}

It is easy to see that $L^2(N_{\mu_0,Q})$ consists of measurable functions $f$ on $\scr P_2$ such that $f\circ\Psi\in L^2(G_Q),$ so that
$$L^2(N_{\mu_0,Q})=\big\{\mu\mapsto  G_Q(\tt f|\Psi)\big|_{\Psi^{-1}(\mu)}:\ \tt f\in L^2(G_Q)\big\},$$
where  $G_Q(\tt f|\Psi)$ is the conditional expectation of $\tt f$   with respect to $G_Q$ given the sigma-algebra $\si(\Psi)$ induced by   $\Psi$, which is constant on the atom $\Psi^{-1}(\mu):=\{h\in T_{\mu_0}: \Psi(h)=\mu\}$ of
$\si(\Psi)$ for each $\mu$.
It is easy to see that $N_{\mu_0,Q}$ is shift-invariant in the following sense.

\beg{prp} Let $\tt h\in T_{\mu_0}$ be a homeomorphism on $\R^d$. Then $N_{\mu_0,Q}= N_{\mu_0\circ\tt h^{-1}, Q\circ \tt h^{-1}}$ for
$Q\circ \tt h^{-1}$ being the linear operator on $T_{\mu_0\circ \tt h^{-1}}$ determined by
$$\big(Q\circ \tt h^{-1}\big)   \tt h_n  := q_n \tt h_n,\ \ n\ge 1,$$
where $\{\tt h_n\}_{n\ge 1}:= \{h_n\circ \tt h^{-1}\}_{n\ge 1}$ is an ONB of $T_{\mu_0\circ \tt h^{-1}}.$ \end{prp}

We have the following result for the OU process on $\scr P_2$.

\beg{thm}\label{T0} Let $(\EE,\D(\EE))$ be defined above. Then the following assertions hold.
\beg{enumerate} \item[ $(1)$] $(\EE,\D(\EE))$ is
  a conservative symmetric Dirichlet form on $L^2(N_{\mu_0,Q})$ with  $\D(\EE)\supset C_b^1(\scr P_2)$ and
\beq\label{A3} \EE(f,g)=\int_{\scr P_2} \<Df(\mu), Dg(\mu)\>_{T_\mu} N_{\mu_0,Q}(\d \mu),\ \ f,g\in C_b^1(\scr P_2).\end{equation}
Moreover, the following log-Sobolev inequality holds:
\beq\label{LS2}  N_{\mu_0,Q}(f^2\log f^2)\le \ff 2 {q_1} \EE(f,f),\ \ f\in \D(\EE),\ N_{\mu_0,Q}(f^2)=1.\end{equation}
\item[$(2)$] The generator $(L,\D(L))$ of $(\EE,\D(\EE))$ has discrete spectrum, and satisfies
\beg{align*}& \D(L)\supset \tt\D(L):= \big\{f\in L^2(N_{\mu_0,Q}):\ f\circ\Psi\in \D(\tt L)\big\},\\
& Lf(\mu) =G_Q\big(\tt L (f\circ\Psi)\big|\Psi=\mu\big):= G_Q\big(\tt L (f\circ\Psi)\big|\Psi\big)\big|_{\Psi=\mu},\ \ f\in \tt\D(L).\end{align*}
\item[$(3)$] Let $P_t$ be the associated Markov semigroup of $(\EE,\D(\EE))$.  Then $P_t$ is compact in $L^2(N_{\mu_0,Q})$ for any $t>0$,
 $P_t$ converges exponentially to $N_{\mu_0,Q}$ in entropy:
\beq\label{EX2}   N_{\mu_0,Q}((P_tf)\log P_t f)\le \e^{-2q_1 t} N_{\mu_0,Q}(f\log f),\ \ t\ge 0,\ 0\le f, \ N_{\mu_0,Q}(f)=1,\end{equation}
and it is hypercontractive:
\beq\label{HP} \beg{split}&\|P_t\|_{L^p(N_{\mu_0,Q})\to L^{p_t}(N_{\mu_0,Q})}:=\sup_{\|f\|_{L^p(N_{\mu_0,Q})}\le 1} \|P_t f\|_{L^{p_t}(N_{\mu_0,Q})}\le 1,\\
&\qquad \ \ t>0, \ p>1, \ p_t:= 1+ (p-1)\e^{2q_1 t}. \end{split}\end{equation}
 \end{enumerate} \end{thm}

\beg{proof} (1) We first prove that $(\EE, \D(\EE))$ is   closed and   \eqref{A3} holds, so that $(\EE, \D(\EE))$ is a symmetric conservative Dirichlet form in $L^2(N_{\mu_0,Q})$.

 Let
$$\EE_1(f):= \EE(f,f)+\|f\|_{L^2(N_{\mu_0,Q})}^2,\ \ \ \tt\EE_1(\tt f) :=\tt\EE(\tt f,\tt f)+\|\tt f\|_{L^2(G_Q)}^2.$$
Let $\{f_n\}_{n\ge 1}\subset \D(\EE)$ such that
$$\lim_{m,n\to\infty} \EE_1(f_n-f_m)= 0.$$ Then
$ f:=\lim_{n\to\infty} f_n$  exists in $L^2(N_{\mu_0,Q})$ and by definition, $\{f_n\circ\Psi\}_{n\ge 1}\subset \D(\tt\EE)$ with
$$\lim_{m,n\to\infty} \tt\EE_1(f_n\circ\Psi-f_m\circ\Psi)= 0.$$ Thus,  the closed property of $(\tt\EE, \D(\tt\EE))$ implies
$$f\circ\Psi=\lim_{n\to\infty} f_n\circ\Psi\in \D(\tt\EE),$$ so that  $f\in \D(\EE)$ by definition. Thus, $(\EE, \D(\EE))$ is   closed.

On the other hand, let $f\in C_b^1(\scr P_2)$. By Proposition \ref{P1} for the reference probability $\mu_0$, we see that for any $h$,
 $$\nn_\phi (f\circ\Psi)(h):=\ff{\d}{\d s}\Big|_{s=0} (f\circ \Psi)(h+s \phi)=  \<\phi, (Df)(\Psi(h))(h)\>_{T_{\mu_0}},\ \ \phi\in {T_{\mu_0}},$$
 so that
 \beq\label{FF}  \nn (f\circ\Psi)(h)= (D f)(\Psi(h))\circ h,\ \ \ \ f\in C_b^1(\scr P_2),\ h\in{T_{\mu_0}}.\end{equation}    Hence, $f\in C_b^1(\scr P_2)$ implies
 $f\circ\Psi\in C_b^1({T_{\mu_0}})\subset \D(\tt \EE),$
 so that  $f\in \D(\EE)$ by definition.

 It remains to verify \eqref{A3}, which together with \eqref{LS1}  implies \eqref{LS2}.   By \eqref{FF} and $N_{\mu_0,Q}= G_Q\circ\Psi^{-1}$,
    \beg{align*}& \EE(f,g):= \tt\EE( f\circ\Psi, g\circ\Psi) =\int_{{T_{\mu_0}}} \<\nn(f\circ\Psi), \nn(g\circ\Psi)\>_{T_{\mu_0}} \d G_Q\\
 &= \int_{T_{\mu_0}} \mu_0\big(\<(D f(\Psi(h)))\circ h, (Dg(\Psi(h)))\circ h\>\big)G_Q(\d h)\\
 &= \int_{T_{\mu_0}} \Psi(h)\big(\<Df(\Psi(h)),Dg(\Psi(h))\>\big) G_Q(\d h)\\
 &= \int_{\scr P_2} \mu\big(\<Df(\mu), Dg(\mu)\>\big) N_{\mu_0,Q}(\d \mu) \\
 &=\int_{\scr P_2} \<Df(\mu), Dg(\mu)\>_{T_\mu} N_{\mu_0,Q}(\d\mu),\ \ \   f,g\in C_b^1(\scr P_2).\end{align*}
  Therefore, \eqref{A3} holds.

(2) As mentioned above that  $(\tt L,\D(\tt L))$ has purely discrete spectrum. So, the set
$$\big\{\tt f\in \D(\tt\EE):\ \tt\EE_1(\tt f):=\tt\EE(\tt f,\tt f)+\|\tt f\|_{L^2(G_Q)}^2\le 1\big\}$$
is compact in $L^2(G_Q)$. By the definitions of $N_{\mu_0,Q}$ and $(\EE,\D(\EE))$, this implies that the set
$$\big\{f\in \D(\EE):\ \EE_1(f):= \EE(f,f)+\|f\|_{L^2(N_{\mu_0,Q})}^2\le 1\big\}$$ is compact in $L^2(N_{\mu_0,Q})$. So,   $L$ has purely discrete spectrum.

Next, let $f\in L^2(N_{\mu_0,Q})$ such that $f\circ\Psi\in \D(\tt L)$. By the definition of $(\EE,\D(\EE))$,  we obtain
\beg{align*}& \int_{\scr P_2} g(\mu) G_Q\big(\tt L(f\circ\Psi)\big|\Psi=\mu\big)N_{\mu_0,Q}(\d\mu)=\int_{T_{\mu_0}} (g\circ\Psi)\tt L(f\circ\Psi)\d G_Q\\
&=-\tt\EE(g\circ\Psi, f\circ\Psi)=-\EE(f,g),\  \  \ g\in\D(\EE).\end{align*}
Thus, $f\in \D(L)$ and $Lf(\mu)= G_Q(\tt L(f\circ\Psi)|\Psi=\mu).$

(3) The log-Sobolev inequality \eqref{LS2}    is equivalent to each of \eqref{EX2} and \eqref{HP},
see \cite{G2}.  Moreover,  by the spectral theory,   since  $L$ has purely discrete spectrum,  $P_t$ is compact  in $L^2(N_{\mu_0,Q})$ for $t>0$.
\end{proof}
\ \newline
{\bf Remark 3.1.} (1) Theorem \ref{T0} implies that $(\EE,C_b^1(\scr P_2))$ is closable in $L^2(N_{\mu_0,Q})$.
We wonder if the closure coincides with $(\EE,\D(\EE))$ or not.

(2)  In general,  $\D(L)\ne  \tt\D(L).$   For any  $f\in \D(L)$, we have
$$ \tt\EE(g\circ \Psi, f\circ\Psi)=\EE(g,f)=-N_{\nu_0,Q}(g (Lf) )= G_Q((g\circ\Psi) (Lf)\circ\Psi),\ \ g\in \D(\EE).$$
If $\{g\circ\Psi: g\in \D(\EE)\}$ is dense in $\D(\tt\EE)$, this would imply that $f\circ\Psi\in \D(\tt L)$ and $\tt L(f\circ\Psi)= (Lf)\circ \Psi$. If so, the OU process on $\scr P_2$ starting at
$\nu_0$ could be constructed as $\nu_t=\Psi(h_t)$ for $h_t$ in \eqref{HT} with $\nu_0=\Psi(h_0)$. However, in general this is not  true, as    there might be
different $h_0$ satisfying $\nu_0=\Psi(h_0)$ and the corresponding  $\Psi(h_t)$ may have different distributions.  In the next section, we formulate $L$ as Laplacian with a drift on $\scr P_2$.

\section{Generator as Laplacian with drift}

  We shall introduce a subclass of $ \D(L)$, such that for functions in this class the generator is formulated as
  $$Lf(\mu)= \DD_{\scr P_2} f(\mu) -  \<b(\mu), D f(\mu)\>_{T_\mu}$$
  for a drift $(b,\D(b))$:  $$b: \scr P_2\supset \D(b) \ni\mu \mapsto b(\mu)\in  T_\mu.$$
  This is compatible with the $d$-dimensional case where the OU process has a generator of type
  $$L_0 f(x) =\DD f(x)- (Ax)\cdot\nn f(x)$$
  for a positive definite $d\times d$-matrix $A$.

\beg{defn}     Let $\D$ be the space of functions $f\in  C_b^2(\scr P_2)$ such that  for $G_Q$-a.e. $h$,
     $$\big(Df(\mu_0\circ h^{-1})\big)\circ h \in \D(Q),\ \ \ \ f\in \D(\DD_{\mu_0\circ h^{-1}}),$$ and
    $$\int_{T_{\mu_0}} \Big|\DD_{\scr P_2} f(\mu_0\circ h^{-1})-\big\<h, Q[(Df(\mu_0\circ h^{-1}))\circ h]\big\>_{T_{\mu_0}} \Big|^2  G_Q(\d h)<\infty. $$     \end{defn}

\paragraph{Remark 4.1.}
(1)    As explained in Remark 2.2(3) that   $\D$  seems too small to be dense    in $L^2(N_{\mu_0,Q})$, so   the following characterization \eqref{AN} on the generator is only formal.
   It would be interesting if one could make this formula meaningful in a weak sense  for a  dense class of functions.

(2) Next, we show     that $\D$ may contain functions given in Remark 2.2(2) of type
   $$f=g(f_1,\dotsc, f_d),\ \ g\in C^\infty_b(\R^d), f_i(\mu):=\mu(\<\cdot,e_i\>).$$
Indeed, by    Remark 2.2(2) we have $f\in \D(\DD_{\scr P_2})$ and $\|\DD_{\scr P_2} f\|_\infty<\infty$. Moreover, noting that
$$f_i(\mu_0\circ h^{-1})=(\mu_0\circ h^{-1})(\<e_i,\cdot\>)=\mu_0(\<e_i,h\>)= \<e_i, h\>_{T_{\mu_0}},$$ we have
\beg{align*} &\big|\big\<h, Q[(Df(\mu_0\circ h^{-1}))\circ h]\big\>_{T_{\mu_0}}\big|^2\\
&= \Big|\sum_{i=1}^d\sum_{m\ge 1}  q_m  \<h,  h_m\>_{T_{\mu_0}} (\pp_i g)(\<h,e_1\>_{T_{\mu_0}},\dotsc, \<h,e_d\>_{T_{\mu_0}}) \<h_m, e_i\>_{T_{\mu_0}}  \Big|^2\\
&\le \bigg(\sum_{i=1}^d \sum_{m\ge 1} \<h,h_m\>_{T_{\mu_0}}^2\bigg)\sum_{i=1}^d \sum_{m\ge 1} q_m^2 (\pp_i g)(\<h,e_1\>_{T_{\mu_0}},\dotsc, \<h,e_d\>_{T_{\mu_0}}) \<e_i, h_m\>_{T_{\mu_0}}^2.\end{align*}
Since
$$\int_{T_{\mu_0}}\bigg(\sum_{i=1}^d \sum_{m\ge 1} \<h,h_m\>_{T_{\mu_0}}^2\bigg)G_Q(\d h)=d\sum_{m\ge 1} q_m^{-1}<\infty,$$
we conclude that $f\in \D$ provided
\beq\label{PP}  \sum_{m\ge 1} q_m^2\<e_i, h_m\>_{T_{\mu_0}}^2<\infty,\ \ 1\le i\le d,\end{equation}
which may be verified by showing that $e_i$ is in a Sobolev type space.
For instance, let $d=1$ (the higher dimension cases can be discussed in the same way) and let $\mu_0 $ be the standard Gaussian measure  such that $L_0:= \DD-x\cdot\nn $ in $L^2(\mu_0)$  has discrete spectrum
$   \{-\ll_{m}:=-  m-1\}_{m\ge 1}  $  and  eigenbasis  $ \{h_{m}\}_{m\ge 1}$ with $ \mu_0(h_m)=0$ for $m\ge 1$.   In this case $e_1=1$ so that $\<e_1,h_m\>_{T_{\mu_0}}=0$ for $m\ge 1$ hence \eqref{PP} holds.

    \beg{thm}\label{G} We have  $  \D\subset \D(L)$ and  for $U_Q^f(h):=  \big\<h, Q[(Df(\mu_0\circ h^{-1}))\circ h]\big\>_{T_{\mu_0}},$
    \beq\label{AN} Lf(\mu) = \DD_{\scr P_2} f(\mu) - G_Q(U_Q^f|\Psi =\mu), \ \ f\in \D. \end{equation} Formally,
    we may write $Lf(\mu)=\DD_{\scr P_2}f(\mu)- \<b(\mu),  Df(\mu)\>_{T_\mu}$ where the drift is given by
    $$ b(\mu)= G_Q\big((Q h)\circ h^{-1}\big|\Psi(h)=\mu\big). $$
\end{thm}

   \beg{proof} (a)   Let $h\in T_{\mu_0}$ and $\mu:=\mu_0\circ h^{-1}$.
  Recall that for any $\tt h\in T_{\mu_0}$,    $\tt h\circ h^{-1}\in T_\mu$ is determined by
 \beq\label{PL} \<\tt h\circ h^{-1}, \phi\>_{T_\mu}= \<\tt h, \phi\circ h\>_{T_{\mu_0}},\ \ \phi\in T_\mu.\end{equation}
Then by Proposition \ref{P1} for the reference probability $\P=\mu_0$, we obtain
\beq\label{A1} \beg{split} &\nn_{\tt h}    (f\circ\Psi)(h) :=\lim_{\vv\downarrow 0} \ff{(f\circ\Psi)(h+\vv\tt h)-(f\circ\Psi)(h)}\vv\\
&= \lim_{\vv\downarrow 0} \ff{f(\mu_0\circ (h+\vv\tt h)^{-1})-f(\mu_0\circ h^{-1})  }\vv=\int_{\R^d} \<Df(\mu)(h), \tt h\>\d\mu_0\\
&= \<Df(\mu), \tt h\circ h^{-1}\>_{T_\mu}= D_{\tt h\circ h^{-1}} f(\mu),\ \ \tt h, h\in T_{\mu_0}, \ \ \mu= \mu_0\circ h^{-1}.\end{split}\end{equation}

(b)  By definition, for any $f\in \D$, $Lf$ given in \eqref{AN} is  a well-defined function in $L^2(N_{\mu_0,Q})$.
    It suffices to prove the integration by parts formula
\beq\label{AN2} \EE(f,g)=-\int_{\scr P_2} g(\mu) L f(\mu) N_{\mu_0,Q}(\d\mu),\ \ g\in C_b^1(\scr P_2).\end{equation}
Simply denote $\tt f=f\circ\Psi$ and $\tt g=g\circ\Psi.$  By \eqref{A3} and the integration by parts formula for $G_Q$, we obtain
\beq\label{AN3}  \beg{split}& \EE(f,g)=\int_{T_{\mu_0}} \<\nn\tt f,\nn\tt g\>\d G_Q=\sum_{n=1}^\infty \int_{T_{\mu_0}} \big(\nn_{h_n}  \tt f\big)(h) \big(\nn_{h_n}\tt g\big)(h) G_Q(\d h)\\
 &= \sum_{n=1}^\infty \int_{T_{\mu_0}} \Big[ \nn_{h_n}  \big( (\nn_{h_n} \tt f)  \tt g \big)(h)-   \tt g(h)\nn_{h_n}\nn_{h_n} \tt f(h) \Big]  G_Q(\d h)\\
 &= \sum_{n=1}^\infty   \int_{T_{\mu_0}}\tt g(h)  \Big[    q_n \<h, h_n\>_{T_{\mu_0}} \big(\nn_{h_n} \tt f\big) (h)     -  \nn_{h_n}\nn_{h_n} \tt f(h) \Big]  G_Q(\d h).\\
 \end{split}\end{equation}
 By \eqref{A1},
\beg{align*} &\sum_{n=1}^\infty q_n \<h, h_n\>_{T_{\mu_0}} \big(\nn_{h_n} \tt f\big) (h) =  \sum_{n=1}^\infty q_n \<h, h_n\>_{T_{\mu_0}} \<h_n, \nn \tt f(h)\>_{T_{\mu_0}}\\
&=  \sum_{n=1}^\infty q_n \<h, h_n\>_{T_{\mu_0}} \<h_n, (D f(\mu_0\circ h^{-1})\circ h)\>_{T_{\mu_0}}
 =\<h, Q [(Df(\mu_0\circ h^{-1}) )\circ h]\>_{T_{\mu_0}},\end{align*}
so that
 \beq\label{AN4}\beg{split} & \sum_{n=1}^\infty   \int_{T_{\mu_0}}\tt g(h)  \    q_n \<h, h_n\>_{T_{\mu_0}} \<h_n, \nn \tt f(h)\>_{T_{\mu_0}} G_Q(\d h)\\
 &= \int_{\scr P_2} g(\mu) G_Q\Big(\big\<h, Q [(Df(\mu_0\circ h^{-1}))\circ h]\big\>_{T_{\mu_0}}\Big|\Psi(h)=\mu\Big) N_{\mu_0,Q}(\d\mu).\end{split} \end{equation}
By Proposition \ref{P1},   \eqref{A1} also implies
   \beq\label{S1}\beg{split}&     \nn_{h_n}\nn_{h_n} \tt f(h)  = \ff{\d}{\d\vv}\Big|_{\vv=0} \mu_0\big(\<(Df)(\mu_0\circ(h+\vv h_n)^{-1})(h+\vv h_n), h_n\>\big)\\
&= \int_{\R^d\times \R^d} \big\<(D^2 f)(\mu_0\circ h^{-1}) (h(x),h(y) ), h_n(x)\otimes h_n(y)\big\>_{HS} \mu_0(\d x)\mu_0(\d y)\\
&\quad + \int_{\R^d}\big\< (\nn Df)(\mu_0\circ h^{-1}) (h(x)),  h_n(x)\otimes h_n(x)\big\>_{HS} \mu_0(\d x)\\
&=I_1(n)+I_2(n),\end{split} \end{equation}
where $\mu=\mu_0\circ h^{-1}.$
\beg{align*} &I_1(n):=  \int_{\R^d\times \R^d} \big\<(D^2 f)(\mu_0\circ h^{-1}) (x,y), (h_n\circ h^{-1})(x)\otimes (h_n\circ h^{-1})(y)\big\>_{HS} \mu (\d x)\mu (\d y),\\
&I_2(n):=\int_{\R^d}\big\< (\nn Df)(\mu_0\circ h^{-1}) (x),  (h_n\circ h^{-1})(x)\otimes (h_n\circ h^{-1})(x)\big\>_{HS} \mu (\d x). \end{align*}
Let $\{\phi_m\}_{m\ge 1}$ be an ONB of $T_\mu:= L^2(\R^d\to \R^d,\mu)$. By \eqref{PL},
$$h_n\circ h^{-1}=\sum_{m\ge 1}  \<h_n\circ h^{-1}, \phi_m\>_{T_\mu}\phi_m =\sum_{m\ge 1}  \<h_n, \phi_m\circ h\>_{T_{\mu_0}}\phi_m,$$
so that
\beg{align*} &\sum_{n=1}^\infty I_1(n)= \sum_{m,l\ge 1}  \sum_{n=1}^\infty \mu_0(\<h_n, \phi_m\circ h\>)\mu_0(\<h_n, \phi_l\circ h\>) \\
&\qquad\qquad\qquad \qquad \cdot \int_{\R^d\times \R^d} \big\<(D^2 f)(\mu) (x,y), \phi_m(x)\otimes \phi_l(y)\big\>_{HS} \mu (\d x)\mu (\d y)\\
&=  \sum_{m,l\ge 1}   \mu_0(\<\phi_m\circ h, \phi_l\circ h\>) \int_{\R^d\times \R^d} \big\<(D^2 f)(\mu) (x,y), \phi_m(x)\otimes \phi_l(y)\big\>_{HS} \mu (\d x)\mu (\d y)\\
&=\sum_{m,l\ge 1}   \mu(\<\phi_m, \phi_l\>) \int_{\R^d\times \R^d} \big\<(D^2 f)(\mu) (x,y), \phi_m(x)\otimes \phi_l(y)\big\>_{HS} \mu (\d x)\mu (\d y)\\
&=\sum_{m\ge 1} \int_{\R^d\times \R^d} \big\<(D^2 f)(\mu) (x,y), \phi_m(x)\otimes \phi_m(y)\big\>_{HS} \mu (\d x)\mu (\d y).\end{align*}
Similarly,
\beg{align*}& \sum_{n=1}^\infty I_2(n)=  \sum_{m,l\ge 1}  \sum_{n=1}^\infty \mu_0(\<h_n, \phi_m\circ h\>)\mu_0(\<h_n, \phi_l\circ h\>)\\
&\qquad\qquad\qquad\qquad \cdot \int_{\R^d}\big\< (\nn Df)(\mu) (x),  \phi_m(x)\otimes \phi_l(x)\big\>_{HS} \mu (\d x)\\
&= \sum_{m\ge 1} \int_{\R^d}\big\< (\nn Df)(\mu) (x),  \phi_m(x)\otimes \phi_m(x)\big\>_{HS} \mu (\d x).\end{align*}
These together with \eqref{S1} and Proposition \ref{LP} yield
 $$\sum_{n=1}^\infty  \nn_{h_n}\nn_{h_n} \tt f(h)= \DD_{\scr P_2} f(\mu_0\circ h^{-1}).$$
 Combining this with \eqref{AN3} and \eqref{AN4}, we prove \eqref{AN}.

 \end{proof}

 \section{Perturbation of the OU process}

 Let $V$ be a measurable function on $\scr P_2$  such that
 $$N_{\mu_0,Q}^V(\d\mu):= \e^{V(\mu)} N_{\mu_0,Q}(\d\mu)$$
 is a probability measure on $\scr P_2$. We consider the pre-Dirichlet form
 $$\EE^V(f,g):=\int_{\scr P_2} \<Df(\mu), D g(\mu)\>_{T_\mu} N_{\mu_0,Q}^V(\d\mu),\ \ f,g\in C_b^1(\scr P_2).$$
 If this form is closable in $L^2(N_{\mu_0,Q}^V)$, then its closure $(\EE^V, \D(\EE^V))$ is a symmetric conservative Dirichlet form, whose generator can be formally written as
 $$L^V f(\mu)= L f(\mu)+\<D V(\mu), Df(\mu)\>_{T_\mu}.$$
 We call the associated Markov process a perturbation of  the OU process.

A simple situation is that $V$ is bounded. In this case,  the closability of $(\EE^V,C_b^1(\scr P_2))$ follows from that of $(\EE, C_b^1(\scr P_2)),$ and \eqref{LS1} implies the log-Sobolev inequality (see \cite{Stroock})
$$N_{\mu_0,Q}^V(f^2\log f^2)\le \ff 2 {q_1} \e^{\sup V-\inf V}\EE^V(f,f),\ \ f\in \D(\EE^V),\ N_{\mu_0,Q}^V(f^2)=1.$$
Consequently, the associate Markov semigroup $P_t^V$ is hypercontractive and exponentially convergent in entropy.  Moreover,
the compactness of $\{f\in \D(\EE): \EE_1(f)\le 1\}$ in $L^2(N_{\mu_0,Q})$ implies that of $\{f\in \D(\EE^V): \EE_1^V(f) \le 1\}$ in $L^2(N_{\mu_0,Q}^V)$, so that the generator $L^V$ has empty essential  spectrum and $P_t^V$ is compact in $L^2(N_{\mu_0,Q}^V)$ for $t>0$. In the following, we intend to extend these to unbounded perturbation $V$.

 In the framework of local Dirichlet forms, unbounded perturbations have been studied in many papers, where
 the key points are  to prove the closability of the pre-Dirichlet form and to see which functional inequalities  of the original Dirichlet form can be kept under the perturbation,
 see for instance \cite{AS, BLW, RZ} and references therein.         However, in all of related references one needs an algebra of bounded measurable functions $\scr A\subset \D(L)$
 which is dense in $\D(L)$  such that the square field is given by
\beq\label{AH} \GG(f,g)= \ff 1 2 \big(L(fg)   -f Lg- gLf\big),\ \ f, g\in \scr A.\end{equation}
In the present situation, the square field reads
$$\GG(f,g)(\mu)=\<Df  (\mu), Dg(\mu)\>_{T_\mu},\ \ \ \ \mu\in \scr P_2,\ f,g\in C_b^1(\scr P_2).$$
But we do not have explicit choice of the  algebra $\scr A$ such that   \eqref{AH} holds.
Therefore, we again come back to  the tangent space $T_{\mu_0}$ by considering the following probability measure on $T_{\mu_0}$:
 $$G_Q^V(\d h):= \e^{(V\circ\Psi)(h)}G_Q(\d h),$$  and the corresponding bilinear form
\beq\label{GPP1} \tt \EE^V(\tt f,\tt g):= \int_{T_{\mu_0}} \<\nn\tt f,\nn\tt g\>_{T_{\mu_0}}\d G_Q^V,\ \ \tt f,\tt g\in C_b^1(T_{\mu_0}).\end{equation} 
 For any $r\in \R$, let $r^+:=\max\{0,r\}$ and $r^-:=(-r)^+.$
By studying  properties of $ \tt\EE^V$,  we obtain the following result under assumption

\beg{enumerate} \item[$(A)$] $V\in C^1(\scr P_2)$ such that $\d N_{\mu_0,Q}^V := \e^V\d N_{\mu_0,Q}$ is a probability measure on $\scr P_2$, and   there exists $p>1$ such that
$$\int_{\scr P_2} \Big(\|DV(\mu)\|_{T_\mu} \e^{V(\mu)^+}+\|D V(\mu)\|_{T_\mu}^p\Big) N_{\mu_0,Q}(\d\mu)<\infty.$$\end{enumerate}

   \beg{thm}\label{TN} Assume $(A)$. Then  the following assertions hold.
   \beg{enumerate}
   \item[$(1)$] $(\EE^V, C_b^1(\scr P_2))$ is closable in $L^2(N_{\mu_0,Q}^V)$, and the closure $(\EE^V,\D(\EE^V))$ is a symmetric conservative Dirichlet form.
  \item[$(2)$] If there exists $\ll>\ff 1{2q_1}$ such that $N_{\mu_0,Q}(\e^{\ll \|DV  \|^2})<\infty$, where $\|DV\|(\mu):=\|DV(\mu)\|_{T_{\mu}},$
  then the associated Markov semigroup $P_t^V$ is
   compact in $L^2(N_{\mu_0,Q}^F)$ for $t>0.$
   \item[$(2)$] If there exists $\vv>0 $ such that
   \beq\label{C1} \int_{\scr P_2} \Big(\e^{\ff{1+\vv}{2q_1} \|DV\|^2 } + \e^{ V^++\vv V^-   } \Big) \d N_{\mu_0,Q} <\infty,\end{equation}
   then there exists a constant $c>0$ such that
   \beq\label{LSN} N_{\mu_0,Q}^V(f^2\log f^2)\le c  \EE^V(f,f),\ \ f\in \D(\EE^V),\ N_{\mu_0,Q}^V(f^2)=1.\end{equation}
   Consequently, for any $t>0$,
   \beg{align*} & \|P_t^V\|_{L^p(N_{\mu_0,Q}^V)\to L^{p_t}(N_{\mu_0,Q}^V)} \le 1,\ \  \  p>1, \ p_t= 1 + (p-1)\e^{4t/c},\\
   &N_{\mu_0,Q}^V\big((P_t^V f)\log P_t^V f\big)\le \e^{-4t/c} N_{\mu_0,Q}^V(f\log f),\ \  \ f\ge 0, \ N_{\mu_0,Q}^V(f)=1.\end{align*}
    \end{enumerate} \end{thm}
Noting that  $N_{\mu_0,Q}=G_Q\circ\Psi^{-1}$,  $\d G_Q^V:= \e^{V\circ\Psi}\d G_Q$  is a probability measure on $T_{\mu_0}$, and
  $$\EE^V(f,g)= \int_{T_{\mu_0}} \<\nn\tt f,\nn\tt g\>_{T_{\mu_0}}\d G_Q^V,\ \ \tt f,\tt g\in C_b^1(T_{\mu_0}),$$
according to the proof of  Theorem \ref{T0}, Theorem \ref{TN}   is a consequence of the following Lemma \ref{LMM} for $F=V\circ\Psi$, where the condition
$ \|\nn F\|_{T_{\mu_0}}\e^{F^+} + \|\nn F\|_{T_{\mu_0}}^p \in L^1(G_Q)$ for some $p>1$ is much weaker than $\e^{\|\nn F\|_{T_{\mu_0}}^2+|F|}\in \cap_{p>1} L^p(G_Q)$
used in \cite[Proposition 3.2]{AS}. We will use the dimension-free Harnack inequality and Bismut formula
for $\tt P_t$ to prove the closability under this weaker condition.  Moreover, the condition \eqref{C1'} for the log-Sobolev inequality is slightly better than that in \cite[Lemma 4.1]{Aida}
where $F^++\vv F^-$ is replaced by $(1+\vv)|F|$.

      \beg{lem}\label{LMM} Let $F\in C^1(T_{\mu_0})$ such that $G_Q^F(\d h):= \e^{F(h)}G_Q(\d h)$ is a probability measure on $T_{\mu_0}$ and
      $ \|\nn F\|_{T_{\mu_0}}\e^{F^+} + \|\nn F\|_{T_{\mu_0}}^p \in L^1(G_Q)$ for some constant $p>1$. Then:
      \beg{enumerate} \item[$(1)$] The bilinear form
   $$\tt\EE^F(\tt f,\tt g):= \int_{T_{\mu_0}} \<\nn \tt f,\nn\tt g\>_{T_{\mu_0}}\d G_Q^F,\ \ \tt f,\tt g\in \F C_b^1(T_{\mu_0})$$
   is closable in $L^2(G_Q^F)$, and the closure $(\tt\EE^F,\D(\tt\EE^F))$ is a symmetric conservative Dirichlet form. Moreover,
   $\D(\tt\EE^F)\supset C_b^1(T_{\mu_0}).$
   \item[$(2)$] If there exists $\ll>\ff 1{2q_1}$ such that $G_Q(\e^{\ll \|\nn F\|_{T_{\mu_0}}^2})<\infty$, then the associated Markov semigroup $\tt P_t^F$ is
   compact in $L^2(G_Q^F)$ for $t>0.$
   \item[$(3)$] If there exists $\vv> 0 $ such that
   \beq\label{C1'} \int_{T_{\mu_0}}  \e^{\ff{1+\vv}{2q_1} \|\nn F\|_{T_{\mu_0}}^2 + F^+ +\vv F^-}   \d G_Q <\infty,\end{equation}
   then there exists a constant $c>0$ such that
   \beq\label{LSN'} G_Q^F  (\tt f^2\log \tt f^2)\le c  \tt \EE^F(\tt f,\tt f),\ \ \tt f\in \D(\tt \EE^F),\ G_Q^F(\tt f^2)=1.\end{equation}
    \end{enumerate}
    \end{lem}

   \beg{proof}  (1) We will establish the integration by parts formula
   \beq\label{ITT} \tt\EE^F(\tt f,\tt g)=-\int_{T_{\mu_0}}  \tt g \big(\tt L^V  \tt f \big)\d G_Q^F,\ \ \tt f,\tt g\in \scr FC_b^2(T_{\mu_0})\end{equation}
   for $\tt L^V\tt f:= \tt L\tt f+\<\nn F,\nn\tt f\>_{T_{\mu_0}},$ so that $(\tt\EE^F, \scr FC_b^2(T_{\mu_0}))$ is closable. Since a function in $\F C_b^1$ can be approximated by functions in $\F C_b^2(T_{\mu_0}) $ under the $C_b^1$-norm,
this also implies that  $(\tt\EE^V, \F C_b^1(T_{\mu_0}))$ is closable.

To this end, we make approximations of $F$. Let $\varphi\in C^\infty(\R)$ such that $\varphi(r)=r$ for $|r|\le 1,$ $1\ge \varphi'\ge 0$ and $\varphi(r)=0$ for $|r|\ge 2.$
For any $m,n\ge 1,$   let
$$F_m:= m\varphi(F/m),\ \ \ F_{m,n}:=\tt P_{\ff 1 n} F_m.$$ We have
\beq\label{IT2} F_m\in C_b(T_{\mu_0})\cap C^1(T_{\mu_0}),\ \ \|\nn F_m\|\le \|\nn F\|.\end{equation}
Since $\tt P_tf(h_0)=\E[f(h_t)]$ for $h_t$ in \eqref{HT}, by \cite[Theorem 3.2.1 and Theorem 3.2.2]{Book} for $A=-Q, b=0$ and $\si(t)=\ss 2$, we have the Harnack inequality
\beq\label{HN} (\tt P_t\tt f(h+v))^p\le (\tt P_t \tt f^p (h ))\e^{\ff{p }{2(p-1)}\|v\|_{T_{\mu_0}}^2},\ \ \tt f\ge 0, p>1, h,v\in T_{\mu_0},\end{equation}
and the Bismut formula
\beq\label{BS} \nn \tt P_t \tt f(h_0)= \ff {\ss 2} t \E\bigg[ \tt f(h_t)\int_0^t \e^{-Qs}  \d W_s\bigg],\ \ t>0, \tt f\in \scr B_b(T_{\mu_0}).\end{equation}
By \eqref{BS}, we see that $F_{m,n}\in C_b^1(T_{\mu_0})$, and   \eqref{HN} together with \eqref{IT2}  and $Q\ge 0$ implies
\beg{align*} &|\tt P_t F_m(h_0+\vv v)- \tt P_t F_m(h_0)|\le \E  |F_m(\vv\e^{-Qt}v+h_t) - F_m(h_t)|\\
&\le \E\int_0^\vv \|\nn F_m\|_{T_{\mu_0}}(r \e^{-Qt}v+h_t)\d r\\
&\le \int_0^\vv (\tt P_t \|\nn F_m\|_{T_{\mu_0}}^p)^{\ff 1 p}(h_0) \e^{\ff{r^2|v|^2}{2(p-1)} }\d r\\
&\le   \vv (\tt P_t \|\nn F \|_{T_{\mu_0}}^p)^{\ff 1 p}(h_0) \e^{\ff{\vv^2|v|^2}{2(p-1)} },\ \ \vv>0, h_0,v\in T_{\mu_0}.\end{align*}
By letting $\vv\downarrow 0$ we derive
\beq\label{HN2} \|\nn F_{m,n}\|_{T_{\mu_0}}=\|\tt P_{\ff 1 n} F_m\|_{T_{\mu_0}}\le (\tt P_{\ff 1 n}\|\nn F \|_{T_{\mu_0}}^p)^{\ff 1 p}.\end{equation}
 Since $F_{m,n} \in C_b^1(T_{\mu_0})$,     the integration by parts formula for    $G_Q$ yields
 \beg{align*} &  \int_{T_{\mu_0}} \<\nn\tt f,\nn \tt g\>_{T_{\mu_0}}\e^{F_{m,n}}  \d G_Q =
 \int_{T_{\mu_0}}\Big(\nn_{\nn\tt f}(\e^{F_{m,n}} \tt g)- \e^{F_{m,n} }\tt g\nn_{\nn \tt f} F_{m,n} \Big)\d G_Q \\
 &= -\int_{T_{\mu_0}} \e^{F_{m,n}} \tt g \big(\tt L \tt f + \nn_{\nn F_{m,n}} \tt f \big)\d G_Q\\
 &= -\int_{T_{\mu_0} }\Big(\tt g(\tt L    + \nn_{\nn F_{m,n}} )\tt f \Big)  \e^{F_{m,n}}\d G_Q,
  \ \ \ \tt f,\tt g\in \F C_b^2(T_{\mu_0}).\end{align*}
  Noting that \eqref{HN2} implies
  $$\Big|\tt g(\tt L    + \nn_{\nn F_{m,n}} )\tt f \Big| \e^{F_{m,n}}\le c_m (1+ \tt P_{\ff 1 n} \|\nn F\|_{T_{\mu_0}}^p)^{\ff 1 p},\ \ n\ge 1  $$
  for some constant $c_m>0$, which are bounded in $L^p(G_Q)$ since
  $$G_Q\big(\tt P_{\ff 1 n} \|\nn F\|_{T_{\mu_0}}^p\big)=G_Q(\|\nn F\|_{T_{\mu_0}}^p)<\infty,$$
  by the dominated convergence theorem we may let  $n\to\infty$ to derive
  $$ \int_{T_{\mu_0}} \<\nn\tt f,\nn \tt g\>_{T_{\mu_0}}\e^{F_{m}}  \d G_Q= -\int_{T_{\mu_0} }\Big(\tt g(\tt L    + \nn_{\nn F_{m}} )\tt f \Big)  \e^{F_m}\d G_Q.$$
  Since \eqref{IT2} and $\e^F+ \|\nn F\|_{T_{\mu_0}}   \e^{F^+}\in L^1(G_Q)$ implies
  $$\big|\tt g(\tt L    + \nn_{\nn F_{m}} )\tt f \big| \e^{F_m} \le c\big(1+\|\nn F\|_{T_{\mu_0}})  \e^{F^+} \in L^1(G_Q),$$
  by using the dominated convergence theorem again, we may let $m\to\infty$ to get \eqref{ITT}.

   Next, let $\tt f\in C_b^1(T_{\mu_0})$, and, for any $n\ge 1$ let
   $$\tt f_n:= \tt f\circ \pi_n,\ \ \pi_n h:= \sum_{i=1}^n \<h,h_i\>_{T_{\mu_0}} h_i.$$
   Then $\{\tt f_n\}_{n\ge 1}\subset \F C_b^1(T_{\mu_0})\subset \D(\tt\EE^F),$ and
   \beg{align*}&\lim_{n\to\infty}  \int_{T_{\mu_0}}\big(|\tt f_n-\tt f|^2+\|\nn(\tt f_n-\tt f)\|_{T_{\mu_0}}^2\big)\d G_Q^F\\
   &= \lim_{n\to\infty} \int_{T_{\mu_0}}\Big(\sum_{l=1}^n \big\{|(\nn_{h_l} \tt f)\circ \pi_n - \nn_{h_l} \tt f|^2 + \sum_{l=n+1}^\infty |\nn_{h_l} \tt f|^2\Big\}  \d G_Q^F\\
   &\le \lim_{n\to\infty}\int_{T_{\mu_0}} \|(\nn  \tt f)\circ \pi_n - \nn  \tt f\|_{T_{\mu_0}}^2 + \sum_{l=n+1}^\infty |\nn_{h_l} \tt f|^2\Big)  \d G_Q^F    =0,\end{align*}
   where the last step follows from  $\nn \tt f\in C_b(T_{\mu_0})$  and the dominated convergence theorem. So, $\tt f \in \D(\tt\EE^F).$

   (2) It suffices to prove that the set
   $$B_1^F:=\big\{\tt f\in C_b^1(T_{\mu_0}):\ \tt \EE_1^F(\tt f):= \tt \EE^F(\tt f,\tt f)+G_Q^F(\tt f^2)\le 1\big\}$$
   is relatively compact in $L^2(G_Q^F)$.
   By the chain rule  and Young's inequality, for any $\vv>0$ and $f\in B_1^F$, we have
   \beg{align*} &\tt\EE_1(\tt f\e^{\ff F 2}, \tt f\e^{\ff F 2})\le G_Q^F(\tt f^2)+ (1+\vv^{-1}) \tt\EE^F(\tt f,\tt f) + \ff{1+\vv} 4 G_Q(\tt f^2 \e^F \|\nn F\|_{T_{\mu_0}}^2)\\
   &\le 1+\vv^{-1} + \ff{1+\vv}{4\ll} G_Q\bigg(\tt f^2\e^F\log\ff{\tt f^2 \e^F }{G_Q(\tt f^2\e^F)}\bigg)+ \ff{1+\vv}{4\ll}G_Q(\tt f^2\e^F)
   \log G_Q\big(\e^{\ll \|\nn F\|_{T_{\mu_0}}^2}\big).\end{align*}
   Since $G_Q(\|\nn(\tt f\e^{\ff F 2})\|_{T_{\mu_0}}^2) \le    \tt\EE_1(\tt f\e^{\ff F 2}, \tt f\e^{\ff F 2})$ and $G_Q(\tt f^2\e^F)= G_Q^F(\tt f^2)\le 1$,
  by combining this with \eqref{LS1} we derive
$$\EE_1(\tt f\e^{\ff F 2}, \tt f\e^{\ff F 2})
    \le 1+\vv^{-1} +\ff{1+\vv }{2q_1\ll}\tt\EE_1(\tt f\e^{\ff F 2}, \tt f\e^{\ff F 2}) +  \ff{1+\vv}{4\ll}\log G_Q\big(\e^{\ll \|\nn F\|_{T_{\mu_0}}^2}\big).$$
  Since $\ll>\ff 1{2q_1}$, we may take a small $\vv>0$ such that $\ff{1+\vv }{2q_1\ll} <1$,  so that this estimate and  $G_Q\big(\e^{\ll \|\nn F\|_{T_{\mu_0}}^2}\big)<\infty$ yield
  $$\EE_1(\tt f\e^{\ff F 2}, \tt f\e^{\ff F 2})\le C,\ \ \ \tt f\in B_1^F$$ for some constant $C>0$.
  Since $\tt L$ has empty essential spectrum, this implies that the set
   $$\{\tt f\e^{\ff F 2}: \tt f\in B_1^F\} $$  is relatively compact in $L^2(G_Q)$. Equivalently,  $B_1^F$  is relatively compact in $L^2(G_Q^F).$

   (3)  The proof of  \eqref{LSN'} is  similar to  that of \cite[Lemma 4.1]{Aida}, but we make a more careful estimate by separating $F^+$ and $F^-$.
   Let $\tt f\in C_b^1(T_{\mu_0})$ such that $G_Q^F(\tt f^2) =1$. By \eqref{LS1} and Young's inequality, we obtain
   \beg{align*} &G^F_Q(\tt f^2\log\tt f^2) =G_Q\big(\tt f^2\e^F \log (\tt f^2 \e^F)\big)- G_Q(\tt f^2 F\e^F)\\
   &\le \ff 2 {q_1} G_Q\big(\|\nn(f\e^{\ff F 2})\|_{T_{\mu_0}}^2\big)+ G_Q(\tt f^2 F^-\e^F)\\
   &\le \ff{2(1+r_1^{-1})}{q_1} \tt\EE^F(\tt f,\tt f) + G_Q^F\bigg(\tt f^2\Big[\ff{1+r_1}{2q_1} \|\nn F\|_{T_{\mu_0}}^2+F^-\Big]\bigg)\\
   &\le \ff{2(1+r_1^{-1})}{q_1} \tt\EE^F(\tt f,\tt f)+ r_2 G_Q^F(\tt f^2\log \tt f^2) + r_2 \log G_Q^F\big(\e^{\ff{1+r_1}{2r_2 q_1} }\|\nn F\|_{T_{\mu_0}}^2+\ff 1 {r_2} F^-\big)\\
   &= \ff{2(1+r_1^{-1})}{q_1} \tt\EE^F(\tt f,\tt f)+ r_2 G_Q^F(\tt f^2\log \tt f^2) + r_2 \log G_Q \Big(\e^{\ff{1+r_1}{2r_2 q_1} \|\nn F\|_{T_{\mu_0}}^2+F^++ \ff {1-r_2} {r_2} F^-}\Big) \end{align*} for any $  r_1,r_2\in (0,1).$
   By taking $r_1 $ small enough and $r_2$ close enough to $r_2$ such that
   $$\ff {1-r_2}{r_2}\lor \Big(\ff{1+r_1}{r_2}-1\Big)  \le \vv,$$  we deduce from this and \eqref{C1'} that
   the defective log-Sobolev inequality
  \beq\label{DFF} G_Q^F(\tt f^2\log\tt f^2)\le c_1 \tt\EE^F(\tt f,\tt f) +c_2,\ \ \tt f\in \D(\tt \EE^F), G_Q^F(\tt f^2)=1 \end{equation}  holds for some constants $c_1,c_2>0$.
By \cite[Theorem  1.1]{AA},  the Dirichlet form $(\tt\EE^F, \D(\tt\EE^F))$ is irreducible.
  According to \cite[Corollary 1.3]{W14} for $\phi(t)=2-t$, \eqref{DFF} implies the log-Sobolev inequality \eqref{LSN'} for some constant $c>0$.
 \end{proof}
  \paragraph{Acknowledgement.}  We would like to thank the referee and Professor L. Dello Schiavo  for very helpful comments and corrections.
 \paragraph{Conflict of interest statement.} The authors have no competing interests as defined by Springer, or other interests that might be perceived to influence the results and/or discussion reported in this paper.
 \paragraph{Data availability statement.}  The results/data/figures in this manuscript have not been published elsewhere, nor are they under consideration (from  one of Contributing Authors) by another publisher.

\end{document}